\documentclass[12pt]{report}
\usepackage[cp1251]{inputenc}
\usepackage[english,russian,ukrainian]{babel}
\usepackage{latexsym,amsfonts,amsthm,amssymb,amsmath}
\usepackage{graphicx,graphics,hhline}
\usepackage{euscript}
\usepackage[all]{xy}
\begin{document}

\Large

\noindent УДК 513.83; 517.5
\vskip 2mm
\smallskip

\noindent{\bf T.M.~Osipchuk (Т.М.~Осіпчук)},  otm82@mail.ru

 \noindent{\bf M.V.~Tkachuk (М.В.~Ткачук)}, mvtkachuk@mail.ru

\vskip 2mm

\noindent Institute of Mathematics NAS of Ukraine

\vskip 4mm

\noindent{\bf ON THE SHADOW PROBLEM

\noindent  AND ITS GENERALIZATIONS TO ELLIPSOIDS}

\vskip 8mm

\small

\noindent The main goal of the paper is to solve some problems about shadow for the sphere generalized on the case of the ellipsoid. Here, the essence of the problem is to find the the minimal number of non-overlapping balls with centers on the sphere which are not holding the center of the sphere and such that  every line passing through the center of the sphere would intersect at least one of the balls. Another method of solving the shadow problem for the sphere is also proposed.

\vskip 4mm

\noindent Основною метою даної статті є  розв'язання  деяких задач про тінь для сфери, узагальнених на випадок еліпсоїда обертання. Тут, суть задачі про тінь полягає в тому, щоб знайти мінімальну кількість куль, що попарно не перетинаються,  з центрами на сфері (кулі не містять її центр) і таких, що довільна пряма, яка проходить через центр сфери, перетинала хоча б одну із вказаних куль. Також, запропоновано інший метод розв'язання задачі про тінь для сфери.

\vskip 4mm

\noindent Главной целью данной статьи является решение некоторых задач о тени для сферы, обобщенных на случай эллипсоида вращения. Здесь, суть задачи о тени состоит в том, чтобы найти минимальное количество попарно непересекающихся шаров с центрами на сфере (и не содержащих ее центр),  таких, чтобы любая прямая, проходящая через центр сферы, пересекала хотя бы один из указанных шаров. Также, предложено другой метод решения задачи о тени для сферы.

\newpage

\large

In 1982 G. Khudaiberganov proposed the problem about shadow \cite{Hud} that can be formulated as follows:  {\it what is the minimal number of closed (open) non-overlapping balls in the space $\mathbb{R}^n$ with radiuses less than the sphere radius,  with centers on it and such that every line passing through the center of the sphere would intersect at least one of the balls.}

This problem (and similar to it) may be written  briefly as: {\it what minimal number of some objects generates the shadow at the point.}

This problem was solved by G. Khudaiberganov for the case $n=2$: it was proved that two balls are sufficient for a circumference on the plane. For all that, it was also made the assumption that for the case $n>2$ the minimal number of such balls is exactly equal to $n$.

Yu. Zelinskii mentioned this problem in many of his works (cf., e.g.  [1\,-\,2]). The problem is also interesting from the standpoint of convex analysis as a special case of the problem about membership of a point to the 1-hull of the union of some collection of balls. In \cite{Zel3}, Yu. Zelinskii and his students proved that three balls are not sufficient for the case $n=3$, but it is possible to generate the shadow at the center of the sphere with four balls. In their work  it is also proved that for the general case the minimal number is $n+1$ balls.   Thus, G. Khudaiberganov's assumption was wrong.  In \cite{Zel3}, it is also proposed another method of solving the problem for the case $n=2$ which gives some numerical estimates. In \cite{Zel4}, the complete answer to this problem for a collection of closed and open balls was obtained.

In the present work, we  provide our own analytically-geometrical  proof of three-dimensional case of the sphere and its generalizations on the prolate ellipsoid. First, we consider the case of the sphere.

Thus, let three balls generate shadow for the center of a sphere. Let us denote them as $B_1$, $B_2$, $B_3$. We also will consider double-napped cones with common vertex at the center of the sphere with generators tangent to the balls. Here and further such cones will be referred to as {\it cones under the balls} and will be denoted as $K_1$, $K_2$, $K_3$, respectively. It is obvious that every line that lies within the cones will cross at least one of the balls. Further, without loss of generality, let us situate the center of the unite sphere at the origin $O$, the center of the firs ball $B_1(A,r_1)$ with the radius $r_1$ at the point $A=(0,0,1)$, $a=\overrightarrow{OA}$, and the center of the second ball $B_2(B,r_2)$ with the radius $r_2$ at the point $B=(0,b_1,b_2)$, $b=\overrightarrow{OB}$.

The conic surfaces $K_1$ and $K_2$ will intersect in two lines which are common tangents to the balls $B_1$, $B_2$ and which will be denoted as $x'$, $x''$ (Figure~1). Those cases when the cones are tangent or do not intersect are not considered for they correspond to such situation of the balls when it is impossible to generate the shadow for the sphere center. Indeed, if cones would not intersect, it could be found a plane passing through the center of the sphere and not intersecting any of the cones. Then the third ball could not overlap the plane. If cones would be tangent then, in the case of open balls, again it could be found the plane passing through the center of the sphere and not intersecting any of the cones and, in the case of closed balls, it would be the same plane without the tangent line which obviously the third ball could not overlap too.

\begin{center}
\includegraphics[width=11 cm]{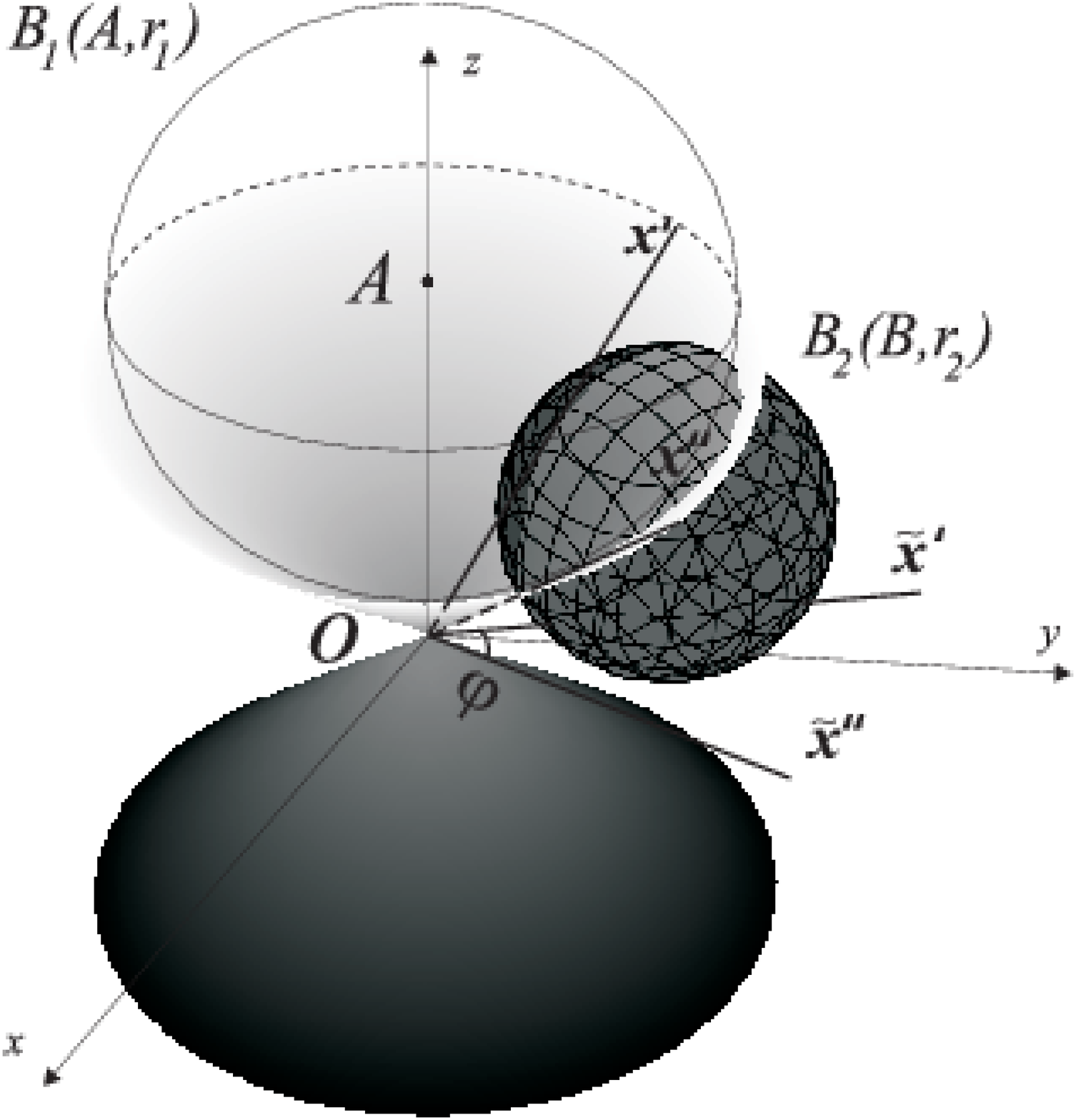}
\end{center}
\small\begin{center} Figure 1 \end{center}
\large

We are interested in the angle $\varphi$ between the projections $\tilde{x}'$, $\tilde{x}''$ of the lines $x'$, $x''$,  on the plane $xOy$, respectively. Let us show that $\angle\varphi<\pi/2$.

Let us consider the point $X=(x_1,x_2,x_3)$ of the intersection of the line $x'$ with the sphere. Then $|x|=1$, where $x=\overrightarrow{OX}$ and

$$
\angle\varphi=2 arctg\frac{x_1}{x_2}.
$$
Let us express $x_1$, $x_2$ in terms of $r_1$, $r_2$.

Using the fact that $|a|=|b|=1$ and $|a-b|=r_1+r_2$ we will get the scalar product
\begin{equation}\label{1}
ab=b_3=1-\frac{1}{2}(r_1+r_2)^2.
\end{equation}

Let the line $x'$ be tangent to the balls $B_1$, $B_2$ at  points $A'$, $B'$, respectively (Figure~2). Then, using the fact that $\cos\widehat{ax}=\sqrt{1-r_1^2}$ and $\cos\widehat{bx}=\sqrt{1-r_2^2}$ obtained from triangles  $OA'A$ and $OB'B$, we get the formulas for the following scalar products:
\begin{equation}\label{2}
ax=x_3=\sqrt{1-r_1^2},
\end{equation}
\begin{equation}\label{3}
bx=b_2x_2+b_3x_3=\sqrt{1-r_2^2}.
\end{equation}
From (\ref{1})--(\ref{3}) and considering $b_2=\sqrt{1-b_3^2}$, we get 
\begin{multline}\label{4}
x_2=\frac{\sqrt{1-r_2^2}-b_3\sqrt{1-r_1^2}}{\sqrt{1-b_3^2}}= \\
=\frac{\sqrt{1-r_2^2}-\left(1-\dfrac{1}{2}(r_1+r_2)^2\right)\sqrt{1-r_1^2}}{\sqrt{(r_1+r_2)^2-\dfrac{1}{4}(r_1+r_2)^4}}.
\end{multline}
And from (\ref{2}) and given $|x|=1$, we obtain
\begin{equation}\label{5}
x_1=\sqrt{1-x_2^2-x_3^2}=\sqrt{r_1^2-x_2^2}.
\end{equation}

\begin{center}
\includegraphics[width=6.5 cm]{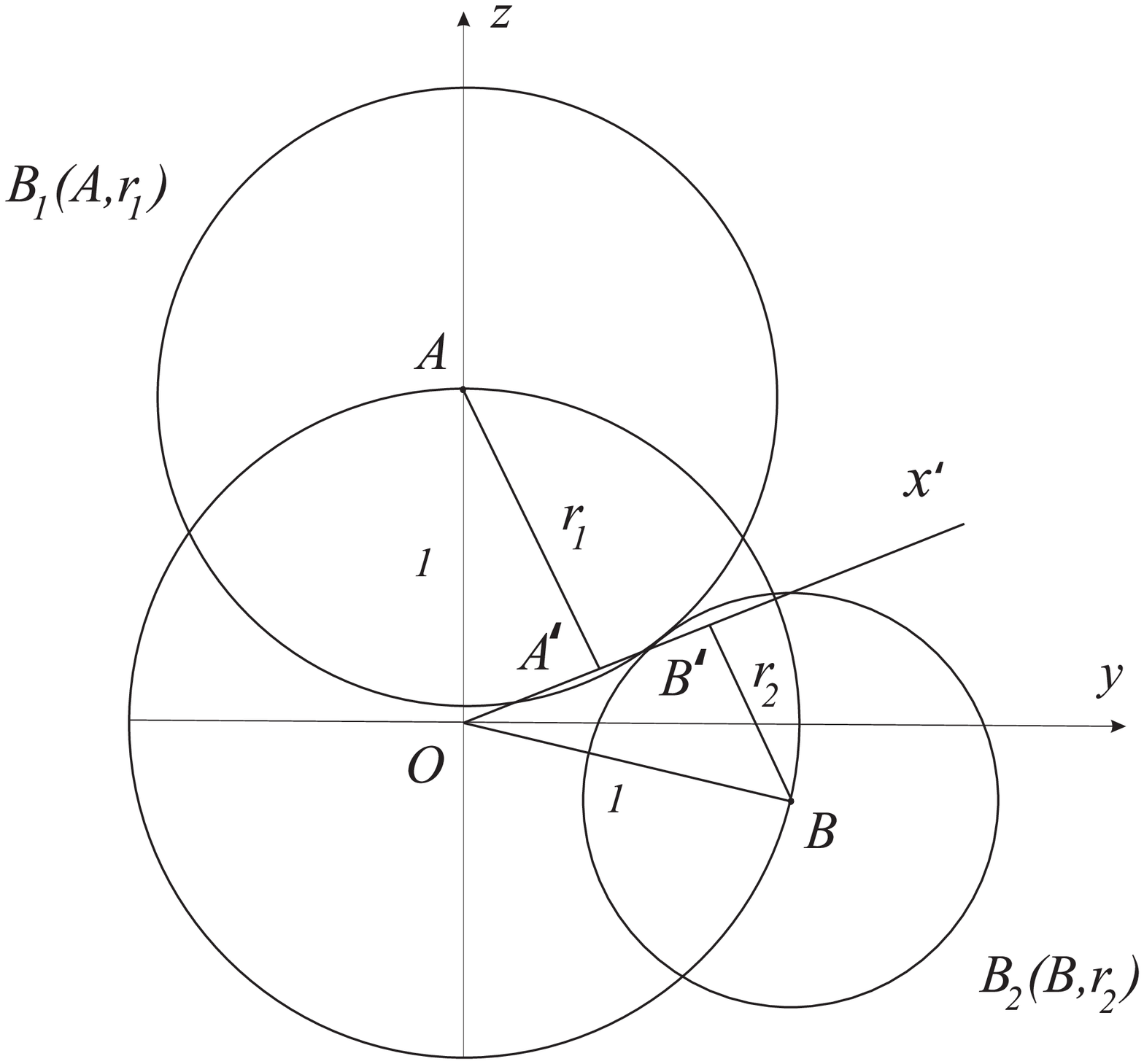}
\end{center}
\small\begin{center} Figure 2 \end{center}
\large

Thus, the expression $\dfrac{x_1}{x_2}$ is a function of two variables $r_1, r_2$ for which a maximum should be found in the open square $(0,1)\,\times\,(0,1)$. Using simple but intensive calculations, we obtain that $\dfrac{x_1}{x_2} < 1$ which implies $\varphi < \pi/2$.

It is obvious that the same considerations are also valid for the ball $B_3$.

Now, let us consider a cone containing the ball $B_1$ and equidistant from the cone $K_1$ on some angle small enough. The whole cone should be contained inside the cones $K_2$, $K_3$, which is impossible, for the projection of the part containing in the cone $K_2$ or  $K_3$ on the plane $xOy$ is in the sector $\varphi < \pi/2$, which gives the whole constructed cone can not be kept in the shadow of the balls $B_2$, $B_3$.

Thus, three balls are not sufficient to generate the shadow for the center of a sphere. But, we remind, that  it can be done with four balls, which is demonstrated in \cite{Zel3}.

Let us consider an interesting fact being  relevant to the problem. Let the balls be open and one of them has the radius that is equal to the radius of the sphere. It is obvious that such a ball generates the shadow everywhere except the equatorial plane tangent to the ball. Thus, for the other balls it would be sufficient to generate the shadow only on the equator. Let the second ball be tangent to the first one, then using the methods of junior course of mathematics,  it is not difficult to show that the maximal angle at which intersection of the second ball with the equatorial plane can be seen from the center of the sphere is obtained in the limit when the center of the second ball tends to the point that is diametrically opposite to the center of the first one and is equal to $\pi/2$. In other words, this fact proves that it is not possible to cover the equator with the other balls.

Now let us consider a prolate ellipsoid instead of the sphere and let us impose the condition such that the open balls with centers on the ellipsoid do not contain the ellipsoid center. If one places the center of the first ball at the base of the ellipsoid minor axis and with the radius equal to this axis, then, similar to the case of sphere, only the equatorial plane stays open. If one chooses ellipsoid major axis long enough and the second ball also situates  at its base, then this ball would cover angle large enough and less then $ \pi$, which gives the possibility for the third ball whose center lies on the minor axis rotation line to cover the rest of the angle.

If, in addition, we consider the closed balls which are not containing the ellipsoid center, then the previous statements
can not be used for the radius of the first ball can not be equal to minor axis. That is why we somewhat decrease it. And, given the possibility for the second ball to be arbitrarily large, we will ensure for the rest of the area to be also closed by the third ball placed on the minor axis rotation line.

 Thus, there exist a prolate ellipsoid and three mutually non-overlapping open (closed) balls with centers on  it which do not intersect the ellipsoid center such that they generate the shadow at this point. It is obvious that two balls are not sufficient  to generate the shadow at the center of any  ellipsoid. Thus, the following proposition is valid.

\noindent {\bf Proposition}. {\it There exists a prolate ellipsoid such that it is necessary and sufficient to have three mutually non-overlapping open (closed) balls with centers on  it which do not intersect the ellipsoid center to generate the shadow at the ellipsoid center.}

Moreover,  for every prolate ellipsoid with ratio between its major and minor axis more than some certain value the Proposition is valid.

The following problem is naturally arising:

\noindent {\bf Problem}. {\it What minimal value of ratio between prolate ellipsoid major and minor axis is sufficient for three mutually non-overlapping open (closed) balls which do not hold the ellipsoid center and with centers on  the ellipsoid to generate the shadow at the ellipsoid center.}

In \cite{Zel2}, there are considered sets similar to the convex ones which are defined as follows.

\noindent {\bf Definition 1.} (\cite{Zel2}) {\it The set  $E\subset\mathbb{R}^n$  is {\it $m$-convex relative to the point }  $x\in \mathbb{R}^n\setminus E$, if there exists an  $m$-dimensional plane $L$ such that $x\in L$  and  $L\cap E=\emptyset$.}

\noindent {\bf Definition 2.} (\cite{Zel2}) {\it The set $E\subset\mathbb{R}^n$  is {\it $m$-convex} if it is $m$-convex relative to every point $x\in \mathbb{R}^n\setminus E$.}

It is easy to verify that both definitions satisfy the axiom of convexity: the intersection of every subfamily of such sets also satisfies the definition. For every set $E\subset\mathbb{R}^n$, it can be considered the minimal $m$-convex set containing $E$  which is called {\it $m$-hull} of the set $E$.

Now, the Problem can be considered as a partial case of the membership of a point to the 1-hull of a union of some collection of balls,  similar to the case of the sphere, as  in \cite{Zel3}, \cite{Zel4}. {\it What minimal value of ratio between major and minor   axis of a prolate ellipsoid ensures the membership of the ellipsoid center to 1-hull of the union of three mutually non-overlapping open (closed) balls which do not intersect the ellipsoid center and with centers on  the ellipsoid.}

We will consider the closed balls for definiteness.

Let us denote radius-vectors of the balls' centers as $a_1, a_2, a_3$ and the balls radiuses as $r_1, r_2, r_3$, respectively. Without loss of generality, let us consider  prolate ellipsoid with the unit minor axis,  the major axis  equal to $d$ and with the center at the origin.

Let $M = \begin{pmatrix} \frac1d & 0 & 0 \\ 0 & 1 & 0 \\ 0 & 0 & 1 \end{pmatrix}$,  then the equation $|Mx| = 1$, where $x$ is the running vector, defines the given ellipsoid.  Under the given conditions,  $|Ma_i| = 1$ and $|a_i - a_j| \geq r_i + r_j$, $i,j = 1,2,3$.

Further, let us build the cones of one nappe under each of the balls with the common vertex at the origin as it was done above. We can get the equation of the cone surface from the formula of cosine of the angle between its axis and generator:
 $$\frac{(a_i, x)}{|a_i||x|} = \frac{\sqrt{|a_i|^2 - r_i^2}}{|a_i|} = \cos \varphi,\quad i=1,2,3.$$
 The equation
 $$\frac{(a_i, x)}{|x|} \geq \sqrt{|a_i|^2 - r_i^2},\quad i\in \{1,2,3\},$$ defines the closed cone.

We consider the system of linear equations
 $$(a_i, x) = \sqrt{|a_i|^2 - r_i^2},\qquad i = 1,2,3.$$
 Let vector $x$ be a solution of the system. If $|x| = 1$, then $x$ lies on the common cone tangent line. If
  $|x| \leq 1$, then the cones under the balls have the common ray $Ox$. If, in addition,  $|x| > 1$, then none of the cones contains $x$, but  $x$ belongs to a common generator of cones  with the smaller cosines of the angles between the cone axis and its generator which contain the cones under the balls in their interior, $$\dfrac{(a_i, x)}{|a_i||x|} = \dfrac{\sqrt{|a_i|^2 - r_i^2}}{|a_i||x|} \leq \dfrac{\sqrt{|a_i|^2 - r_i^2}}{|a_i|},$$   therefore, the intersection of the cones under the balls contains only the point $O$.

  Now, we also consider three cones symmetrical to the given ones with respect to the origin.  For the 1-hull of three balls to contain the origin, it is necessary and sufficient that the union of given six cones contain the whole space. For that, what is obvious, it is sufficient that each three cones which do not have the pare of symmetrical  cones posses nontrivial intersection.

 Let us denote the collection of eight vectors as  $\{s_i = (\pm 1, \pm 1, \pm 1)\}$. Then all conditions of the respective three cones' intersections can be presented as  $(a_i, x_k) = s_k \sqrt{|a_i|^2 - r_i^2}$, $|x_k| \leq 1, k = 1,...,8, i=1,2,3$.

We obtain the following problem of the conditional minimization of function:


$(a_i, x_k) = s_k \sqrt{|a_i|^2 - r_i^2}$, $|x_k| \leq 1$,

$s_k = (\pm 1, \pm 1, \pm 1)$, $k = 1,...,8$, $i = 1,2,3$,

$|Ma_i| = 1$, $|a_i - a_j| \geq r_i + r_j$, $r_i \geq 0$, $1 \leq i < j \leq 3$,

$M = \begin{pmatrix} \frac1d & 0 & 0 \\ 0 & 1 & 0 \\ 0 & 0 & 1 \end{pmatrix}$,

$d \to \min.$

Lagrange multiplier method or numerical methods can be applied to this problem.

Since all restrictions are algebraic, we may use Sequential Least SQuares Programming from the package SciPy for Python to find the minimum. The calculations give us the minimum value approximately equal to $2\sqrt2$.

Minimum $d_{min} = 2 \sqrt2$ can not be reached and it is the limit value for the following three-parameter family of configurations:

$r_1=x, r_2=y, a_1=(0,0,1), a_2=(0,u,v),$

$u=-1+\frac12 (4-(x+y)^2), v = \frac12 (x+y)\sqrt{4-(x+y)^2},$

$a_3(z,0,0), r_3=\sqrt{z^2+1}-1,$

when  $x \to 1, y \to 1, z \to 2\sqrt2$, where $r_1$, $r_2$, $r_3$ are radiuses of the balls, and $x$,~$y$,~$z$~are the parameters.

Therefore,  the desired minimal  value of ratio between the  ellipsoid major and minor axis equals to $2 \sqrt2$.

Thus, the following theorem is proved:

\noindent {\bf Theorem 1.} {\it Let a prolate ellipsoid with value of ratio between its major and minor axis more then $ 2 \sqrt2$ is given. In order that the center of the ellipsoid belongs to $1$-hull of the union of mutually non-overlapping closed (open)  balls which do not hold the ellipsoid center and with centers on  the ellipsoid, it necessary and sufficient to have three balls.}

Let three balls generate the shadow at the center of a prolate ellipsoid. Let us change the scale along the ellipsoid major axis to transform, by virtue of linear conversion, the ellipsoid to the sphere and the balls to the prolate ellipsoids with parallel axis. Now we use the invariance of  $1$-convexity with respect to affine transformation (\cite{Zel2}), which gives the obtained ellipsoids to generate the shadow at the center of the sphere.

In such a way,  the following theorem is proved:

\noindent {\bf Theorem 2.} {\it   It is sufficient to have three mutually non-overlapping closed (open)   prolate ellipsoids which do not hold the sphere center and with centers on  the sphere to generate the shadow at the sphere center.  Moreover, the ellipsoids can be mutually homothetic.}

\renewcommand{\bibname}{References}

\end{document}